\documentclass{amsart}
\usepackage{amsmath}

\usepackage{mypackage_old}
\usepackage{hyperref}
\usepackage{tikz}

\newtheorem{theorem}{Theorem}
\newtheorem{corollary}{Corollary}

\newtheorem{definition}{Definition}

\newtheorem{lemma}{Lemma}
\newtheorem{proposition}{Proposition}

\hypersetup{colorlinks=true,
		  linkcolor=blue}

\newcommand{\Ha}{\mathcal H}
\newcommand{\La}{\mathscr L}

\newcommand{\clop}{\mathrm{clopen}}

\newcommand{\Neg}{\mathcal N}
\newcommand{\circo}
{\,\raisebox{1pt}{\tikz \draw[line width=0.5pt] circle(1.2pt);}\,}

\begin{document}

\mybic
\date \today 
 
\title
[The vector integral]
{Some Properties of The Finitely Additive Vector Integral}

\subjclass{
Primary 28B05, 28C05; Secondary 28C15, 46A55}
\keywords	{
Bochner integral,
Choquet integral representation,
Finitely additive probabilities,
Integral representation,
Pettis integral,
Radon measure.
}

\begin{abstract} 
We prove some results concerning the finitely additive, 
vector integrals of Bochner and Pettis and their representation
over a countably additive probability space. An application
to the non compact Choquet theorem is also provided.
\end{abstract}

\maketitle

\section{Introduction}

In this paper we are concerned with various instances of 
the condition
\begin{equation}
\label{intro}
\int_\Omega h(f)dm
	=
\int_S h(\tilde f)d\tilde m,
\qquad
h\in\Ha
\end{equation}
involving two finitely additive probabilities, $m$ and 
$\tilde m$, two functions $f$ and $\tilde f$ with values 
in a Banach space $X$ and an appropriate family $\Ha$ 
of real valued functions on $X$ (the case of functions 
$h\in\Ha$ with values in another Banach space may be 
treated similarly but does not seem to add much value). 
In the typical situation considered in this work, the pair 
$(m,f)$ and the family $\Ha$ are taken as given and the 
problem is finding a pair $(\tilde m,\tilde f)$ which solves 
\eqref{intro}, a {\it representation} of  $(m,f)$ relatively 
to $\Ha$ (see Definition \ref{def rep}). 

The problem of whether $\tilde f$ may represent $(m,f)$
for some probability $\tilde m$ and relatively to a sufficiently
interesting family $\Ha$ of functions is solved in section
\ref{sec cover} under fairly general conditions involving
the range of $f$ and of $\tilde f$. Our interest, however,
goes beyond mere existence and aims at representations
satisfying some additional, natural properties. In particular 
we focus on conditions hinging on the supporting set $S$ 
and on the intervening measure $\tilde m$.

In section \ref{sec N} we consider representations 
supported by the set $\N$ of natural numbers. We 
show (Theorem \ref{th bochner}) that this kind of 
representation is always possible when $f$ is 
measurable. If, in addition, $f$ is integrable then,
passing from $\N$ to its compactification $\beta\N$,
we deduce a countably additive representing measure.

In section \ref{sec abstract} we examine the problem 
of the existence of a countably additive representing 
measure $\tilde m$ in the case in which $f$ is just 
Pettis integrable. This case is at the same time more 
interesting and more delicate. In fact, outside of some
special cases, e.g. when $X$ is reflexive or a dual space
or when $f$ is tight, there is no natural representation 
$(\tilde m,\tilde f)$. Stone space techniques are only 
partly useful  and, in particular, do not permit to identify 
$\tilde f$.

From a historical prospective, the first result concerning 
countably additive representations of finitely additive 
integrals was proved by Dubins and Savage 
\cite{dubins_savage} in the setting $X=S=\R=\Omega$ 
and with $\Ha$ consisting of bounded, continuous 
functions on the real line. In their representation 
$f=\tilde f$. The most relevant result so far (of  which 
I am aware) was obtained by Karandikar 
\cite{karandikar_88} (but see \cite{karandikar} 
as well) with the purpose of extending limit theorems
to the finitely additive setting. In these works $X=\R^\N$ 
and $\mathcal H$ is formed by all bounded continuous 
functions which only  depend on finitely many coordinates. 
In this framework, a key observation is that measurable 
functions with values in $\R^n$ are tight, which makes 
it possible to apply the theorems of Dini and of Daniell. 
We mention that, motivated again by limit theorems 
applications, this result was extended by Berti, Regazzini
and Rigo \cite[Theorem 2.1]{berti_et_al_07} to the case 
of functions with values in a metric space but assuming 
tightness explicitly. In infinite dimensional spaces, in 
fact, measurable functions need not be tight nor need 
continuous transformations of measurable functions 
be measurable.

In section \ref{sec RNP} we consider the case in which 
$X$ possesses the Radon-Nikod\`ym property ({\it RNP}) 
and $f$ is Pettis integrable. In Theorem \ref{th pettis} we 
prove the existence of a countably additive representation 
for functions on $X$ which are uniformly continuous with
respect to the weak topology. In Theorem \ref{th pettis choice}
we find conditions, hinging on the range of $f$, under 
which a countably additive representation obtains with
$\tilde f=f$. In section \ref{sec PIP} we discuss the Pettis
integrability property ({\it PIP}), rarely discussed in the
context of finite additivity, obtaining some partial conclusions,
particularly under the assumption that $X$ is separable.
Eventually, in section \ref{sec choquet} we prove a version
of Choquet integral representation for non compact, non
convex sets. This final result permits to appreciate 
advantages and disadvantages of our approach.

\section{Preliminaries}

Throughout the paper $(\Omega,\A,m)$ will be a fixed
finitely additive probability space, with $\Omega$ an 
arbitrary non empty set, $\A$ an algebra of subsets of 
$\Omega$ and $m$ a finitely additive probability defined 
on $\A$, i.e. $m\in\Prob\A$ (the symbol $\Prob[ca]\A$ 
denotes countably additive probabilities)%
\footnote{
We write $\Prob S$ for probabilities defined on the power
set of a set $S$.}. 
Also given 
will be a Banach space $X$, its dual space $X^*$ (with 
$\mathbb B$ and $\mathbb B^*$ the closed unit balls 
in $X$ and $X^*$ respectively) and a function $f$ of 
$\Omega$ into $X$ on which, from time to time, 
different measurability and integrability properties 
will be imposed (relatively to $m$). 

Convergence in measure and measurability are defined 
differently under finite additivity and it is useful to recall 
these definitions explicitly: a sequence $\seqn f$ of maps 
from $\Omega$ to $X$ converges in $m$ measure to 
$0$ (or simply $m$-converges to $0$) whenever
\begin{equation}
\label{converge}
\lim_nm^*(\norm{f_n}>c)=0,
\qquad
c>0.
\end{equation}
In \eqref{converge} $m^*$ is the outer measure induced 
by $m$, i.e.
\begin{equation}
m^*(E)
	=
\inf\{m(A):A\in\A,\ E\subset A\}.
\end{equation}
$f$ is measurable, in symbols $f\in L^0_X(m)$ (or 
$f\in L^0_X(\A,m)$ if reference to $\A$ is not obvious 
from the context), if there exists a sequence $f_n$ of 
$\A$ simple, $X$ valued functions whose distance 
from $f$ $m$-converges to 0. When $m$ is countably 
additive and $\A$ a $\sigma$ algebra, this notion of 
measurability implies the classical $(\A(m),\Bor(X))$ 
definition (with $\A(m)$ denoting the $m$ completion 
of $\A$ and $\Bor(X)$ the $\sigma$ algebra of Borel 
subsets of $X$) and that the two are actually equivalent 
if the range of $f$ is essentially separably valued, 
\cite[III.6.10]{bible}. We speak of $f$ as integrable, 
and write $f\in L^1_X(m)$, whenever there exists a 
sequence $f_n$ of $\A$ simple functions which 
$m$-converges to $f$ and such that 
$\lim_{m,n}\int\norm{f_m-f_n}dm=0$%
\footnote{
Some authors speak of functions measurable or integrable
in the above defined sense as strongly measurable and
Bochner integrable. We will not follow this terminology.
}. 
Measurability and integrability of $f$ are qualified as 
{\it weak} or {\it norm} when the corresponding property 
characterizes the set $\{x^*f:x^*\in\mathbb B^*\}$ or
the function $\norm f$, respectively. 

Although not strictly necessary, it will spare a considerable 
amount of repetitions to assume, as we shall do with no 
further mention, that $f$ is weakly and norm measurable. 

The definite Pettis integral of $f$ over $A\in\A$, if it exists, is 
the unique element $x_A\in X$ satisfying 
$x^*x_A=\int_Ax^*fdm$ 
for all $x^*\in X^*$. Then $f$ is Pettis integrable if it 
admits a definite Pettis integral over every $A\in\A$.

Concerning function spaces, the symbol $\Fun{S,T}$ 
(resp. $\Fun S$) denotes the functions from $S$ to $T$ 
(resp. from $S$ to $\R$) and if $A\subset S$ and 
$g\in\Fun{S,T}$, the image of $A$ under $g$ is 
indicated by $g[A]$. If $S$ is a topological space, 
$\cl F$ denotes the closure of $F\subset S$, while
$\Cts S$ (resp. $\Cts[u]S$) indicates the family of
real valued, continuous (resp. uniformly continuous) 
functions on $S$. The suffix $b$ attached to a class 
of real valued function indicates the subclass consisting 
of bounded functions (e.g. $\Cts[ub]S$). We shall 
use repeatedly the following version of a result of 
Hager \cite[Theorem 4.2]{hager}: given a family 
$\mathscr G\subset\Fun{X}$ of linear fuctionals,
the collection of compositions $U(g_1,\ldots,g_n)$ 
with $g_1,\ldots,g_n\in\mathscr G$ and 
$U\in\Fun[u]{\R^n}$ is uniformly dense in
$\Cts[u]{(X,\tau_{\mathscr G})}$ where 
$\tau_{\mathscr G}$ is the initial topology induced 
by $\mathscr G$.

We will often make reference to the normed linear spaces
\begin{equation}
\label{La}
\La(m,f)
	=
\left\{H\in\Fun{\Omega\times X}:
H(\cdot,f(\cdot))\in L^1(m)\right\},
\end{equation}
with  
$\norm[\La(m,f)]H
	=
\norm[L^1(m)]{H(f)}$%
\footnote{
Abusing notation we shall identify $\Fun X$ with the 
corresponding subspace of $\Fun{\Omega\times X}$.
We shall also write for simplicity $H(f)$ in place of
$H(\cdot,f(\cdot))$.
},
and $\La_u(m,f)=\Cts[u]X\cap\La(m,f)$.

\begin{definition}
\label{def rep}
Let $S$ be a non empty set and $\Ha\subset\La(m,f)$. 
Then $(\tilde m,\tilde f)$ is a representation of $(m,f)$ 
relatively to $\Ha$ and with support $S$ if 
$\tilde f\in\Fun{S,X}$, $\tilde m\in\Prob S$,
$\Ha\subset\La(\tilde m,\tilde f)$ and
\begin{equation}
\int h(f)dm
=
\int h(\tilde f)d\tilde m,
\qquad
h\in\Ha.
\end{equation}
\end{definition}

\section{A general theorem on representations}
\label{sec cover}

We now prove a general result on representations which 
will be useful in the sequel. It gives an abstract characterization 
of this kind of problems and it is  perhaps of its own interest. 
Given that this result only depends on uniform continuity 
and that both the strong and the weak topology make $X$ 
into a uniformizable space, in the context of this section 
we will consider $X$ only for its properties as a uniform 
space. We refer to \cite{willard} for terminology on uniform 
spaces.

\begin{definition}
\label{def cover}
Let $\mathscr D$ be a diagonal uniformity on $X$. A set 
$K\subset X$ is an $m$-cover of (the range of) $f$ (or 
$(m,\mathscr D)$-cover) if
\begin{equation}
\sup_{D\in\mathscr D} m^*\big(f\notin D(K)\big)=0
\end{equation}
where $D(K)=\bigcup_{k\in K}\{x:(x,k)\in D\}$. 
We speak of an $X$ valued function $g$ as a $m$-cover
of $f$ if the range of $g$ $m$-covers $f$.
\end{definition}

If $X$ is metric, then the sets $D(K)$ will just be the
cover of $K$ with open balls of given radius. In terms
of Definition \ref{def cover} a $m$-null function is 
simply an $X$ valued function which is $m$-covered 
by the origin.

\begin{theorem}
\label{th cover}
Let $\mathscr D$ be a diagonal uniformity on $X$ 
and $\Ha=\Cts[u]{(X,\mathscr D)}\cap\La(m,f)$. Let 
$S$ be a non empty set, $\Neg$ an ideal of its subsets 
and $\tilde f\in\Fun{S,X}$. Consider the following 
properties:
\begin{enumerate}[(a)]
\item
\label{pro cover}
$\tilde f[N^c]$ is an $(m,\mathscr D)$-cover of $f$ 
for each $N\in\Neg$;
\item
\label{pro represent}
there exists $\tilde m\in\Prob{S}$ vanishing on 
$\Neg$ such that $(\tilde m,\tilde f)$ represents 
$(m,f)$ relatively to $\Ha$.
\end{enumerate}
Then \imply{pro cover}{pro represent} and, if $\mathscr D$ 
is the weak uniformity generated by $\Ha$, 
\imply{pro represent}{pro cover}.
\end{theorem}

\begin{proof}
Assume \iref{pro cover} and fix $h\in\Ha$. For each
$N\in\Neg$, $D\in\mathscr D$ and $n\in\N$ there 
exists $A_n\in\A$ such that $m(A_n)>1-2^{-n}$ and 
$f[A_n]\subset D(\tilde f[N^c])$. Thus,
\begin{align*}
\int h(f)dm
&=
\lim_n\int_{A_n}h(f)dm
\\&\ge
\inf_{z\in D(\tilde f[N^c])}h(z)
\\&\ge
\inf_{s\in N^c}h(\tilde f)(s)-\sup_{x,z\in D}\abs{h(z)-h(x)}.
\end{align*}
Given that $h\in\Cts[u]{(X,\mathscr D)}$ we conclude 
\begin{equation*}
\sup_{N\in\Neg}\inf_{s\in N^c}h(\tilde f)(s)
\le
\int h(f)dm.
\end{equation*}
Claim \iref{pro represent} then follows from 
\cite[Theorem 4.5]{JCA_2018}.

Assume now that $(\tilde m,\tilde f)$ is as in 
\iref{pro represent}. Fix $N\in\Neg$ and let 
$K=\cl{\tilde f[N^c]}$.
Of course, to check that $\tilde f[N^c]$ is an
$(m,\mathscr D)$-cover of $f$ it is enough to
restrict attention to a base of the uniformity
$\mathscr D$. If $\mathscr D$ is the weak uniformity 
induced by $\Ha$, one such base consists of sets of the
form
\begin{equation}
D_0
=
\bigcap_{i=1}^n
\big\{(x,y)\in X\times X:\abs{h_i(x)-h_i(y)}<\varepsilon_i\big\},
\qquad
h_1,\ldots,h_n\in\Ha.
\end{equation}
Then, $y\notin D_0(\tilde f[N^c])$ if and only if 
\begin{equation}
1
	=
\inf_{s\in N^c}\sup_{\{i=1,\ldots,n\}}
\dabs{\frac{h_i\big(\tilde f(s)\big)-h_i(y)}{\varepsilon_i}}
\wedge 1
	=
H(y).
\end{equation}
Notice that $H\in\Ha$, $H=0$ on $\tilde f[N^c]$ and 
$H=1$ on $D_0(\tilde f[N^c])^c$. But then
\begin{align*}
0
&=
\int_{N^c}H(\tilde f)d\mu 
=
\int H(\tilde f)d\mu
=
\int H(f)dm
\ge
m^*\big(f\notin D_0(\tilde f[N^c])\big).
\end{align*}
\end{proof}

It is implicit in Theorem \ref{th cover} that obtaining 
representations by embedding $\Omega$ or  the range 
of $f$ in some larger set is a trivial exercise which gains 
some interest only if it permits to obtain special properties 
on the representing measure $\tilde m$, such as countably 
additivity. This is the case of compactifications that we 
shall consider in the following sections.

Although Theorem \ref{th cover} is very useful in proving 
existence of representations, it is also too general to 
provide detailed information on the representing pair. 
This will then require additional structure.

In applications of Theorem \ref{th cover} we will
omit reference to $\mathscr D$ whenever this is the
uniformity induced by the norm.

\begin{corollary}
\label{cor cover}
Let $f\in L^0_X(m)$ and fix a non empty set $S$. 
Then $(m,f)$ admits a representation relative to 
$\La_u(m,f)$ and supported by $S$ if and only if 
$f$ has an $m$-cover of cardinality not greater 
than $S$.
\end{corollary}

\begin{proof}
Let $K\subset X$ be an $m$-cover of $f$. If the cardinality 
of $K$ is not greater than that of $S$ then there exists 
$\tilde f\in\Fun{S,X}$ which is onto $K$. It is then clear 
from Theorem \ref{th cover}	that $(\tilde m,\tilde f)$ 
is the claimed representation of $(m,f)$ for some 
$\tilde m\in\Prob S$. For the converse observe that if
$f\in L^0_X(m)$ then $\Cts[ub]X\subset\La(m,f)$ and 
that the norm uniformity coincides with the weak 
uniformity generated by $\La_u(m,f)$
(see e.g. \cite[37.9]{willard}).
\end{proof}

Corollary \ref{cor cover} clarifies the role of measurability
in this class of problems. If $f$ is measurable then all 
uniformly continuous transformations of $f$ are 
measurable too and the set $\La_u(m,f)$ is sufficiently 
rich to generate the norm uniformity. With weak 
measurability we can only establish representations 
relative to the family $\Cts[u]{(X,weak)}\cap \La(m,f)$. 
The diagonal uniformity induced by this class of functions 
is the same as that induced by $X^*$. The sets $D(K)$ 
defining an $m$-cover would then just be the weakly 
open sets. See Corollary \ref{cor bochner Gamma}
below.

\section{Representations over $\N$.}
\label{sec N}

By Theorem \ref{th cover} if
$X$ is separable, then $(m,f)$ admits a representation
supported by $\N$. More interestingly, if $\seqn f$ is 
a sequence of $\A$ simple functions $m$-converging 
to $f$, then $\bigcup_nf_n[\Omega]$ is an $m$-cover 
of $f$. Thus, measurable functions admit representations
supported by $\N$. 

The prominence of $\N$ as a prototypical model for
{\it finitely additive} integration was clearly noted, 
among others, by Maharam \cite{maharam_76}. It 
appears from the preceding remark that it is a model 
for vector integration {\it in general}. Among other 
features of $\N$, the classical decomposition of 
Yosida and Hewitt takes an especially simple form 
since the set of countably additive measures on 
$\N$ is isomorphic to $\ell^1$ while a purely finitely 
additive set function is characterized by the property 
of vanishing on singletons (i.e. of being weightless, 
in Maharam's terminology which we adopt). An 
example of weightless, atomless measures are limit 
frequencies, which are often suggested as a convenient 
statistical model and are the most appropriate model
for the uniform distribution over the integers, see e.g. 
Kadane and O'Hagan \cite{kadane_ohagan}. 

In the following Theorem \ref{th bochner} we get, 
further to existence, a more precise description of 
the representation. Remarkably, we find that the 
representing probability may always be chosen 
to be weightless. This finding seems to contrast with 
the popular view which considers this family of set 
functions virtually useless.

\begin{theorem}
\label{th bochner}
Let $f\in L^0_X(m)$. There exist 
$\kappa\in\Fun{\N,\Omega}$, a countable algebra 
$\A_0\subset\A$ and $\mu\in\Prob{\N}$ weightless
such that, letting $F=f\circo\kappa$, the following is
true: $f\in L^0_X(\A_0,m)$, 
$F\in L^0_X(\kappa^{-1}(\A_0),\mu)$,
$\La_u(m.f)\subset\La(\mu,F)$ and
\begin{equation}
\label{bochner}
\int_A h(f)dm
	=
\int_{\kappa^{-1}(A)} h(F)d\mu,
\quad
A\in\A_0,\ h\in\La_u(m,f)
\end{equation}
Moreover, if $f$ is Bochner, norm or weakly integrable 
then so is $F$.
\end{theorem}

\begin{proof}
Fix a sequence $\seqn g$ of $\A$ simple functions such 
that
$m^*(\norm{f-g_n}>2^{-n})
	\le
2^{-(n+1)}$
and choose $D_n\in\A$ such that 
$D_n\subset\{\norm{f-g_n}\le2^{-n}\}$
and
$m(D_n)>1-2^{-n}$.
Fix the trivial partition $\pi_0$ and define inductively 
$\pi_n$ as the partition induced by $\pi_{n-1}$, $D_n$
and the sets supporting $g_n$. Denote by $\A_n$ the 
algebra generated by $\pi_n$ and $\A_0=\bigcup_n\A_n$. 
For each $n\in\N$ and $E\in\pi_n$ choose 
$\omega_E^n\in E$ 
and define 
\begin{equation}
\label{f simple}
f_n
	=
\sum_{E\in\pi_n}f(\omega^n_E)\set E,
\qquad
n\in\N.
\end{equation}
Of course,
$\norm{f-f_n}
	\le
2^{-(n-1)}$
on $D_n$. Thus $\seqn{f}$ is a sequence of $\A_0$ 
simple functions $m$-converging to $f$ and 
$f\in L^0_X(\A_0,m)$. 

Let $\kappa\in\Fun{\N,\Omega}$ be an enumeration of
$\{\omega^n_E:E\in\pi_n,n\in\N\}$
and set $F=f\circo\kappa$. Enumerate $\A_0$ as 
$\{A_1,A_2,\ldots\}$. Fix $k\in\N$. For all $j\in\N$ 
sufficiently large and all $A_1,\ldots,A_i\in\A_j$ then
\begin{equation}
\label{inclusion}
(f_j,\set{A_1},\ldots,\set{A_i})[\Omega]
	\subset
\big(F,\set{\kappa^{-1}(A_1)},\ldots,\set{\kappa^{-1}(A_i)}\big)
[\{k,k+1,\ldots\}].
\end{equation} 
Consider the $X\oplus c_0$ valued functions 
$\hat f$ and $\hat F$ implicitly defined as
$\hat f(\omega)
=
(f(\omega),\ldots,2^{-i}\set{A_i}(\omega),\ldots)$
and
$\hat F(n)
=
(F(n),\ldots,2^{-i}\set{\kappa^{-1}(A_i)}(n),\ldots)$.
From the inequality
\begin{align*}
\bnorm[X\oplus c_0]{\hat f(\omega)-\hat F(n)}
&=
\norm[X]{f(\omega)-F(n)}
+
\sup_i2^{-i}\babs{\set{A_i}(\omega)-\set{\kappa^{-1}(A_i)}(n)}
\\&\le
\norm[X]{f_j(\omega)-F(n)}
+
\sup_{i\le I}2^{-i}
\babs{\set{A_i}(\omega)-\set{\kappa^{-1}(A_i)}(n)}
\\&\quad+
2^{-I}
+
\norm[X]{f_j(\omega)-f(\omega)}
\end{align*}
we conclude
\begin{equation}
\inf_{n\ge k}
\bnorm[X\oplus c_0]{\hat f(\omega)-\hat F(n)}
	\le
2^{-I}
+
\norm[X]{f_j(\omega)-f(\omega)}.
\end{equation}
Consequently, $\hat F[\{k,k+1,\ldots\}]$ $m$-covers 
$\hat f$ for each $k\in\N$ and, by Theorem \ref{th cover}, 
there exists $\mu\in\Prob S$ such that $(\mu,\hat F)$ 
represents $(m,\hat f)$ relatively to 
$\La(m,\hat f)\cap\Cts[u]X$ and that 
$\mu(\{1,\ldots,k\})=0$ for any $k\in\N$. 
Thus $\mu$ is weightless. Notice that if $A\in\A_0$
and $h\in\Ha$ then $h(f)\set A=\hat h(\hat f)$ for 
some $\hat h\in\La(m,\hat f)\cap\Cts[u]{X\oplus c_0}$.
From this we deduce \eqref{bochner}. The claim
that $F\in L^0_X(\kappa^{-1}(\A_0),\mu)$ is clear
from the first part of the proof. 

Concerning the last claim, given that $f$ and $F$ 
are measurable, Bochner, norm or weak integrability 
depend on which functions are included in $\Ha$. 
The corresponding property thus carry over from 
$f$ to $F$.
\end{proof}

The only explicit link between $m$ and $\mu$ is 
the relation $\mu(\kappa^{-1}(A))=m(A)$ for each 
$A\in\A_0$. Thus, in case $m$ is countably additive, 
then so will be the restriction $\mu_0$ of $\mu$ to 
$\kappa^{-1}(\A_0)$. However, by \cite[Corollary 3]
{ascherl_lehn}, $\mu_0$ has then a countably additive 
extension $\mu_1$ to the power set  of $\N$ and 
such extension would definitely satisfy \eqref{bochner}. 
We have an example in which a countably additive 
and a purely finitely additive set function induce 
exactly the same representation. 

\begin{corollary}
Let $(S,\Sigma,\lambda)$ be a countably additive
probability space. Then $L^1_X(\lambda)$ is
isometrically isomorphic with a closed subset of 
$\ell^1_X$.
\end{corollary}

\begin{proof}
Choose $g\in L^1_X(\lambda)$. By Theorem \ref{th bochner} 
and the preceding remark there exists $\mu_1\in\Prob[ca]\N$ 
and a representation $G$ of $g$ supported by $\N$ such that, 
\begin{equation}
\int h(g)d\lambda
=
\sum_nh(G)(n)\mu_1(n),
\qquad
h\in\Cts[u]X\cap\La(\lambda,g).
\end{equation}
Since $g$ is Bochner integrable, so is $G$. Let
$T(g)_n=G(n)\mu_1(n)$, for any $n\in\N$.
Clearly, 
$\norm[\ell^1_X]{T(g)}
=
\sum_n\norm{G(n)}\mu_1(n)
=
\int\norm gd\lambda
=
\norm[L^1_X(\lambda)]g$.
If $\seqn g$ is a sequence in $L^1_X(\lambda)$
such that $T(g^n)$ converges in $\ell^1_X$ then
$\seqn g$ is Cauchy in $L^1_X(\lambda)$. Its
norm limit $g_0$ is such that $T(g_0)$ is the
norm limit of $T(g_n)$.
\end{proof}

Given the completely general nature of the Banach 
space $X$, Theorem \ref{th bochner} admits some 
easy extensions to the case in which $X$ is the direct
sum of a finite or a countable family of Banach spaces.
These extensions imply that the representation obtained
in Theorem \ref{th bochner} may be constructed so as
to preserve $m$-convergence or convergence in $L^1(m)$.

\begin{corollary}
\label{cor bochner}
For each $i\in\N$, let $X_i$ be a Banach space and 
$f_i\in L^0_{X_i}(m)$. There exist 
$\kappa\in\Fun{\N,\Omega}$, a countable algebra
$\A_0\subset\A$ and $\mu\in\Prob\N$ weightless 
such that $f_i\in L_X^0(\A_0,m)$ and, letting 
$F_i=f_i\circo\kappa$,
\begin{equation}
\label{bochner multi}
\int_A h(f_1,\ldots,f_n)dm
	=
\int_{\kappa^{-1}(A)} h(F_1,\ldots,F_n)d\mu
\end{equation}
for each $n\in\N$, $A\in\A_0$ and $h\in
\bCts[u]{\oplus_{i=1}^nX_i}
\cap
\La(m,\oplus_{i=1}^nf_i)$,
\end{corollary}

\begin{proof}
Fix $N\in\N$. Of course $Y_N=\oplus_{i=1}^NX_i$ is a 
Banach space if endowed, e.g., with the norm 
$\norm{x_1,\ldots,x_N}
	=
\sup_{1\le i\le N}\norm[X_i]{x_i}$.
Moreover, the function $\oplus_{i=1}^Nf_i$ belongs to
$L^0_{Y_N}(m)$. We can apply Theorem \ref{th bochner} 
and find a map $\kappa_N\in\Fun{\N,\Omega}$ such that 
the range of $(\oplus_{i=1}^Nf_i)\circo\kappa_N$ in restriction 
to all cofinite sets $m$-almost covers $(\oplus_{i=1}^Nf_i)$. 
This conclusion remains valid if we replace $\kappa_N$ 
with the enumeration $\kappa$ of 
$\bigcup_{N\in\N}\kappa_N[\N]$. Thus we obtain the same 
change of variable for each $N$. The claim follows noting 
that, if $F_i=f_i\circo\kappa$, then 
$\oplus_{i=1}^NF_i
=
(\oplus_{i=1}^Nf_i)\circo\kappa$.
\end{proof}

Thus, if $X=X_1=X_2=\ldots$ and $f_n$ converges to $f_1$ 
either in measure or in $L_X^1(m)$ then $F_n$ converges to 
$F_1$ in the corresponding topology. The extension to an
infinite direct sum requires to specify the norm explicitly.
We adopt the following definition:
\begin{equation}
\bigoplus_{i\in\N}X_i
=
\left\{x\in\bigtimes_{i\in\N}X_i:
\sum_i2^{-i}\norm[X_i]{x_i}<\infty\right\}.
\end{equation}

\begin{corollary}
\label{cor bochner multi}
For each $i\in\N$ let $X_i$ be a Banach space 
and $f_i\in L^1_{X_i}(m)$. Assume that 
$\sup_i\norm[X_i]{f_i}\in L^1(m)$. 
Then there exist $\kappa\in\Fun{\N,\Omega}$ 
and a countable subalgebra $\A_0\subset\A$ 
such that, letting $F_i$ and $\mu$ be as in 
Corollary \ref{cor bochner}, $F_i\in L^1_{X_i}(\mu)$ 
and
\begin{equation}
\label{bochner multi}
\int_A h(f_1,f_2,\ldots)dm
	=
\int_{\kappa^{-1}(A)} h(F_1,F_2,\ldots)d\mu,
\end{equation}
for every $A\in\A_0$ and 
$h\in\Cts[u]{\oplus_iX_i}\cap\La(m,\oplus_if_i)$.
\end{corollary}

\begin{proof}
Put $b=\sup_i\norm[X_i]{f_i}$. There is no loss of generality 
in assuming $\int bdm=1$ and $\{b=+\infty\}=\emp$. We claim 
that $f=\oplus_if_i$ is an element of $L^1_X(m)$. 

Indeed,
$
\norm[\oplus_iX_i]{f(\omega)}
	=
\sum_i2^{-i}\norm[X_i]{f_i(\omega)}
	\le 
b(\omega)
$
so that $f\in\Fun{\Omega,X}$.
For each $i\in\N$, let $\sseqn{g_i^n}$ be a sequence
of $\A$ simple functions such that
$\int\norm[X_i]{f_i-g_i^n}dm
	\le
2^{-n}$. 
If we  let $\pi^n$ be the meet of the finite, $\A$ measurable
partitions associated with $g_i^n$ for $i=1,\ldots,n$ then 
we can write
\begin{equation}
\label{gn}
g_i^n
	=
\sum_{D\in\pi^n} \chi_i^n(D)\set{D},
\qtext{with}
\chi_i^n\in\Fun{\pi^n,X_i}
\qquad
i=1,\ldots,n
\end{equation}
or else $\chi_i^n=g_i^n=0$ if $i>n$. The function
$g^n
	=
\oplus_ig_i^n$
maps $\Omega$ into $X$ and is $\A$ simple.
In fact if 
$\chi^n
	=
\oplus_i\chi_i^n
$
we may represent $g^n$ as
\begin{equation}
g^n
	=
\sum_{D\in\pi^n}\chi^n(D)\set{D}.
\end{equation}
Then,
$\norm[\oplus_iX_i]{f-g^n}
	=
\sum_{i\le n}2^{-i}\norm[X_i]{f_i-g_i^n}
+
\sum_{i>n}2^{-i}\norm[X_i]{f_i}
	\le
\sum_{i\le n}2^{-i}\norm[X_i]{f_i-g_i^n}
+
2^{-n}b$, 
so that 
\begin{equation}
\int\norm[\oplus_iX_i]{f-g^n}dm\le2^{-(n-1)}.
\end{equation}
The claim then follows from Theorem \ref{th bochner}.
\end{proof}

The condition $f\in L^0_X(m)$ turns out to be quite
restrictive in applications. For example, the identity 
map on a subset $A\subset X$ is measurable if and 
only if for each $\varepsilon>0$ there is a finite 
collection of open balls of radius $\varepsilon$ such 
that the part of $A$ which cannot be covered by 
such family has probability less then $\varepsilon$. 
A case in which measurability may be disposed of is 
given next.

\begin{corollary}
\label{cor bochner Gamma}
Let $f$ be norm integrable and let $X^*$ be separable. 
Then Theorem \ref{th bochner} applies upon replacing 
$\La_u(m,f)$ with $\Cts[u]{(X,weak)}\cap\La(m,f)$.
\end{corollary}

\begin{proof}
Let $\{x^*_n:n\in\N\}$ be dense in $\mathbb B^*$. Define 
$\Lambda\in\Fun{X,\ell^1}$ by letting
\begin{equation}
\label{Lambda}
\Lambda(x)_n
	=
x^*_n(x)2^{-n},
\qquad
n\in\N.
\end{equation}
Given that $f$ is weakly and norm measurable, by 
our general assumption, then 
$\Lambda(f)\in L^0_{\ell^1}(m)$. 
It follows from Theorem \ref{th bochner} that there 
exists $\kappa\in\Fun{\N,\Omega}$ and $\A_0$ 
such that, letting $F=f\circo\kappa$, one gets
\begin{equation}
j(\Lambda(F))\in L^0_X(\mu)
\qand
\int_Aj(\Lambda(f))dm
	=
\int_{\kappa^{-1}(A)}j(\Lambda(F))d\mu,
\end{equation}
for every $A\in\A_0$ and $j\in\Cts[ub]{\ell^1}$. The 
compositions $j(\Lambda)$ include all functions of 
the form $\alpha(x^*_{n_1},\ldots,x^*_{n_k})$, with 
$x^*_{n_1},\ldots,x^*_{n_k}\in\mathbb B^*$ and 
$\alpha\in\Cts[ub]{\R^{k}}$ and are thus dense in
$\Cts[ub]{(X,weak)}$ with respect to the uniform 
topology. The extension to
$\Cts[u]{(X,weak)}\cap\La(m,f)$
is obvious. 
\end{proof}

An alternative proof may be obtained by noting 
that if $X^*$ is separable then the weak topology 
of $\mathbb B$ is metrizable and the proof of 
Theorem \ref{th bochner} carries over with only 
minor modifications.

In anticipation of the following sections, a standard
application of Stone space techniques delivers a
countably additive representation.

\begin{theorem}
\label{th bochner stone}
Assume that $f\in L^1_X(m)$.
Let $\chi$ be the Stone isomorphism of $\A$ and
the field $\F$ of clopen sets of its Stone space $S$. 
Let $\lambda\in\Prob[ca]{\sigma\F}$ be the
extension of $m\circo\chi^{-1}$ . There exists 
$\tilde f\in L^1_X(\lambda)$ such that 
$\La_u(m,f)
\subset 
\La_u(\lambda,\tilde f)$
and that, for each $A\in\A$,
$\lambda(\chi(A)\cap\{\tilde f\notin\cl{f[A]}\})=0$
and 
\begin{equation}
\label{bochner stone}
\int_Ah(f)dm
	=
\int_{\chi(A)}h(\tilde f)d\lambda,
\qquad
h\in\La_u(m,f).
\end{equation}
If $f$ has closed range there exists 
$\nu\in\Prob[ca]{\sigma\{h(f):h\in\Cts[u]{X}\}}$ 
such that 
\begin{equation}
\label{bochner stone daniell}
h(f)\in L^1(\nu)
\qand
\int h(f)dm
	=
\int h(f)d\nu,
\qquad
h\in\La_u(m,f).
\end{equation}
\end{theorem}

\begin{proof}
Clearly, the isomorphism $\chi$ extends to $\A$ simple
functions by letting
\begin{equation}
\chi\Big(\sum_{i=1}^nx_i\set{A_i}\Big)
	=
\sum_{i=1}^nx_i\set{\chi(A_i)}.
\end{equation}
If $\seqn g$ is a sequence of $\A$ simple functions
that converges in $L^1_X(m)$ to the origin then 
$\chi(g_n)$ forms a Cauchy sequence in 
$L^1_X(\lambda)$ and its limit cannot be but  the 
origin. Thus we obtain a further extension of $\chi$
as a map of $L^1_X(m)$ into $L_X^1(\lambda)$. Write
$\tilde f=\chi(f)$. Let $\seqn f$ be a sequence of $\A$ 
simple functions converging to $f$ in $L^1_X(m)$. Fix 
$A\in\A$, $\delta>0$ set 
$C_n=\{\norm{\tilde f-\chi(f_n)}<\delta\}$
and choose $E_n\in\A$ such that 
$E_n\subset\{\norm{f-f_n}<\delta\}$.
Then
\begin{align*}
\tilde f[C_n\cap\chi(A\cap E_n)]
&\subset
\chi(f_n)[\chi(A\cap E_n)]+\delta\mathbb B
\\&=
f_n[A\cap E_n]+\delta\mathbb B
\\&\subset
f[A]+2\delta\mathbb B.
\end{align*}
Thus, 
$\lambda^*(\chi(A)
\cap\{\tilde f\notin f[A]+2\delta\mathbb B\})
\le
\lambda(C_n^c)+m(E_n^c)$ 
and we conclude that
$\lambda\big(\chi(A)\cap\big\{\tilde f\notin\cl{ f[A]}\big\}\big)=0$.

Fix $h\in\La_u(m,f)$ and $A\in\A$. Then,
$h(f_n)$ converges to $h(f)$ in $L^1(m)$ and $h(\chi(f_n))$ 
converges to $h(\tilde f)$ in $L^1(\lambda)$. But then,
\begin{align*}
\int_A h(f)dm
	&=
\lim_n\int_{A}h(f_n)dm
	\\&=
\lim_n\int_{\chi(A)}h(\chi(f_n))d\lambda
	\\&=
\int_{\chi(A)}h(\tilde f)d\lambda.
\end{align*}

If $f$ has closed range and $\seqn h$ is a sequence in 
$\La_u(m,f)$ such that $h_n(f)$ decreases 
pointwise to $0$, then $h_n$ decreases pointwise to 
$0$ on $f[\Omega]$ i.e. $h_n(\tilde f)$ decreases to 
$0$ $\lambda$ a.s.. It follows from \eqref{bochner stone} 
that $\int h(f)dm$ is a Daniell integral over the vector 
lattice $\La_u(m,f)$ which contains the unit. 
\end{proof}

Theorem \ref{th bochner stone} requires the assumption 
$f\in L^1_X(m)$ which is quite restrictive and is hardly 
satisfied in applications. The reason is that Stone 
isomorphism would otherwise not be sufficient to
identify a representation $\tilde f$. A similar problem
will be encountered in the next section. On the other 
hand if $f$ is just Pettis integrable Theorem 
\ref{th bochner stone} is no longer valid.

\section{Abstract standard representations.}
\label{sec abstract}

In the present and the following section we look for
a countably additive representation without assuming
that $f$ is measurable and by exploiting the Stone-\v{C}ech 
compactification. It turns out that this technique is less 
useful than one may expect. In fact the compactification 
$\beta\Omega$ of $\Omega$ does not permit to find a 
natural extension of $f$. In particular it is hardly possible 
to obtain a continuous extension that may serve 
as a representation with respect to some countably 
additive measure. This is however possible in the 
following special case

\begin{theorem}
\label{th dual}
Let $X$ be the dual space of some Banach space $W$
and assume that $f$ is norm integrable. Then there 
exists a Radon probability measure $\lambda$ defined 
on $\Bor(\beta\Omega)$ and a map 
$\tilde f\in\Fun{\beta\Omega,X}$ such that
\begin{equation}
\label{dual}
\int h(f)dm
=
\int h(\tilde f)d\lambda,
\qquad
h\in\Cts[u]{(X,weak^*)}\cap\La(m,f).
\end{equation}
\end{theorem}

\begin{proof}
Consider the family $\Fun[b]{\Omega,X}$ of functions 
whose range is contained in some multiple of $\mathbb B$ 
(the unit sphere of $X$). Endowing $X$ with the weak$^*$ 
topology and $\Omega$ with the discrete topology, every 
$g\in\Fun[b]{\Omega,X}$ is continuous with values 
in a compact, Hausdorff space. It admits a continuous 
extension 
$\tilde g
	\in
\bFun{\beta\Omega,X}$
(see e.g. \cite[19.5]{willard}). There exists a Radon 
probability measure $\lambda$ defined on 
$\Bor(\beta\Omega)$ such that
\begin{equation}
\int h(g)dm
	=
\int h(\tilde g)d\lambda
\end{equation}
for all $g\in\Fun[b]{\Omega,X}$ and all 
$h\in\Cts{(X,weak^*)}\cap\La(m,g)$.
Define
\begin{equation}
f_n
=
\frac{f}{1+2^{-n}\norm f}.
\end{equation}
Clearly, $f_n\in\Fun[b]{\Omega,X}$. Moreover,
$\norm{f-f_n}$ converges to $0$ in $L^1(m)$.
Then, letting $A$ run across finite subsets of 
the unit sphere of $W$, and taking advantage of
$\tau$ additivity of $\lambda$
\begin{align*}
0
&=
\lim_{m,n}\sup_A
\int\sup_{w\in A}\abs{(f_n-f_m)w}dm
\\&=
\lim_{m,n}\sup_A
\int\sup_{w\in A}\abs{(\tilde f_n-\tilde f_m)w}d\lambda
\\&=
\lim_{m,n}\int\norm{\tilde f_n-\tilde f_m}d\lambda.
\end{align*}
We may find a subsequence (still indexed by $n$ for
convenience) such that 
\begin{align*}
\infty
>
\lim_k
\int\sum_{n=1}^k\norm{\tilde f_{n+1}-\tilde f_n}d\lambda
=
\int\sum_n\norm{\tilde f_{n+1}-\tilde f_n}d\lambda.
\end{align*}
There exists therefore a $\lambda$ null set outside of 
which $\tilde f_n$ is convergent in $X$ to some limit
$\tilde f$ such that $\norm{\tilde f}\in L^1(\lambda)$
and that $\norm{\tilde f_n-\tilde f}$ converges to $0$
in $L^1(\lambda)$. It is also clear that $\tilde f$ is
independent of the chosen subsequence. If 
$U\in\Cts[ub]{\R^p}$ and $w_1,\ldots,w_p\in W$
\begin{align*}
\int U(w_1,\ldots,w_p)(f)dm
&=
\lim_n\int U(w_1,\ldots,w_p)(f_n)dm
\\&=
\lim_n\int U(w_1,\ldots,w_p)(\tilde f_n)d\lambda
\\&=
\int U(w_1,\ldots,w_p)(\tilde f)d\lambda.
\end{align*}
Therefore \eqref{dual} follows from Hager's lemma.
\end{proof}

In the general case, however, the situation is more 
complicated. We found it convenient to compactify 
$\Omega\times X$. This has the double advantage 
of finding a natural extension $T(H)$ of the elements 
$H$ of $\La(m,f)$ and, on the other hand, to arrive 
at the countably additive representation of the 
integrals of such extensions. We refer to this result as 
an {\it abstract standard representation} of $(m,f)$, 
see \eqref{T int}. Its proof is given in the following 
Lemma \ref{lemma T}.  

Nevertheless, the problem of expressing $T(H)$ in 
the form $H(\tilde f)$ for some $X$ valued function 
$\tilde f$ proves to be quite difficult to solve. The
more so given our choice to treat $f$ just a weakly
and strong measurable function. To this end we
later consider several additional assumptions on the
space $X$.

The following is essentially an application of Stone%
-\v{C}ech and Stone-Weierstrass to our setting. We
need to spell out several properties.

\begin{lemma}
\label{lemma T}
There exist: 
(i)
a countably additive probability space $(S,\Bor(S),\lambda)$ 
with $S$ compact and Hausdorff and $\lambda$ regular,
(ii)  
a Boolean homomorphism $\chi$ of $\A(m)$ 
into $\clop(S)$ and
(iii)
an isometric, vector lattice homomorphism $T$ of 
$\La(m,f)$ into $L^1(\lambda)$ such that
\begin{equation}
\label{T int}
\int_A H(f)dm
	=
\int_{\chi(A)} T(H)d\lambda,
\qquad
A\in\A,\ 
H\in\La(m,f).
\end{equation}
In addition: 
\begin{enumerate}[(a).]
\item
\label{onto C}
$T$ is an algebraic isomorphism of $\La_b(m,f)$ onto $\Cts S$;
\item\label{meas}
each $\tilde f\in\Cts S$ is $\sigma(\clop(S))$ measurable;

\item
\label{clop}
if $E\in\clop(S)$ then $\lambda(E\bigtriangleup\chi(A))=0$
for some $A\in\A(m)$;

\item
\label{dense}
$T[\La(m,f)]$ is norm dense in $L^1(\lambda)$;

\item
\label{sw}
for arbitrary $n\in\N$, $H_1,\ldots,H_n\in\La(m,f)$ and 
$\alpha\in\Cts[b]{\R^n}$
\begin{equation}
\label{T commute}
\alpha\big(T(H_1),\ldots,T(H_n)\big)
	=
T\big(\alpha(H_1,\ldots,H_n)\big),
\qquad
\lambda\ a.s.;
\end{equation}

\item\label{range}
for each $A\in\A(m)$ and $h\in\La(m,f)$, $h(f)$ is 
$m$-null on $A$ if and only if $T(h)$ is $\lambda$-null
on $\chi(A)$; moreover, 
$\lambda(\chi(A)\setminus\tilde A)=0$ where
\begin{equation}
\label{tilde A}
\tilde A
	=
\chi(A)\cap\bigcap_{\{h(f)\text{ $m$-null on $A$}\}}\{T(h)=0\}.
\end{equation}
\end{enumerate}
\end{lemma}

\begin{proof}
Since $\La_b(m,f)$ is a closed subalgebra of $\Fun[b]{\Omega\times X}$ 
containing 
the unit as well as a vector lattice there exist \cite[IV.6.20]
{bible} a compact, Hausdorff space $S$ and an isometry 
$U$ between $\La_b(m,f)$ (endowed with the uniform norm) 
and $\Cts S$. The map $U$ is also an isomorphism of 
algebras as well as of partially ordered sets. Define 
\begin{equation}
\chi(A)
	=
\{U(\set{A\times X})=1\},
\qquad
A\in\A(m).
\end{equation}
Of course, $0\le U(\set A)=U(\set A)^2\le 1$ so  that
$U(\set A)
	=
\set{\chi(A)}$
and $\chi(A)$ is clopen. 

The integral $\int H(f)dm$ for $H\in\La_b(m,f)$ 
may thus be rewritten as $\phi(UH)$ with $\phi$
a linear functional on $\Cts S$ with 
$\norm\phi=\norm m=1$.  By Riesz-Markov there exists
$\lambda\in\Prob[ca]{\Bor(S)}$ regular and such that
for any $A\in\A$ and any $H\in\La_b(m,f)$
\begin{equation}
\label{U}
\int_A H(f)dm
	=
\int U(\set AH)d\lambda	
	=
\int_{\chi(A)} U(H)d\lambda.
\end{equation}
The extension $T$ of $U$ to $\La(m,f)$ is obtained by 
first letting $T_0(H)=\lim_k U(H\wedge k)$ for $H\ge0$, 
then noting that $T_0$ is additive on $\La(m,f)_+$ 
(because
$(G_1+G_2)\wedge k
	\le 
G_1\wedge k+G_2\wedge k
	\le
(G_1+G_2)\wedge 2k$)
and eventually defining, up to a $\lambda$ null set,
\begin{equation}
\label{T}
T(G)
	=
T_0(G^+)-T_0(G^-),
\qquad
G\in\La.
\end{equation}
Clearly,  $T$ is a linear map of $\La(m,f)$ into 
$L^1(\lambda)$ and its restriction to $\La_b(m,f)$ 
coincides with $U$ and is therefore an isomorphism 
of algebras. This proves \iref{onto C}. To prove the 
other properties of $T$, using linearity and the fact 
that if $G_1,G_2\in\La(m,f)$ have disjoint supporting 
sets, then so have $T(G_1)$ and $T(G_2)$, it is enough 
to restrict to $\La(m,f)_+$. Thus, for any choice of 
$G_1,G_2\in\La(m,f)_+$ we see that
\begin{align*}
\int_AG_1(f)dm
	&=
\lim_k\int_A\big[G_1(f)\wedge k\big] dm
	\\&=
\lim_k\int_{\chi(A)}U(G_1\wedge k)d\lambda
	\\&=
\int_{\chi(A)}T(G_1)d\lambda
\end{align*}
(so that \eqref{T int} holds) and that
\begin{equation*}
T(G_1\wedge G_2)
	=
\lim_kU(G_1\wedge G_2\wedge k)
	=
\lim_kU(G_1\wedge k)\wedge U(G_2\wedge k)
	=
T(G_1)\wedge T(G_2).
\end{equation*}
It is clear from \eqref{T int} that 
$\norm[L^1(\lambda)]{T(H)}
	=
\norm[\La(m,f)]H$.

\iref{meas}
Let $\tilde f\in\Cts S$ and $H\in\La_b(m,f)$ be such that 
$T(H)=\tilde f$. Given that $H(f)$ is the $L^1(m)$ 
limit of a sequence of $\A$ simple functions, then $\tilde f$
is the $L^1(\lambda)$ limit of a sequence of $\chi(\A)$
simple functions.

\iref{clop}.
Let $E\in\clop(S)$. Define $E_1=\{U^{-1}(\set E)=1\}$
and $E_2=\{\omega:(\omega,f(\omega))\in E_1\}$.
Since $\set E\in\Cts S$, then by \iref{onto C} 
$\set{E_1}\in\La_b(m,f)$. Moreover, 
$\set{E_1}(\omega,f(\omega))
=
\set{E_2}(\omega)$
so that $E_2\in\A(m)$ and 
 $\lambda(E\bigtriangleup\chi(E_2))=0$.

\iref{dense}.
The claim follows from the inclusion
$L^\infty(\lambda)
	\subset
\cl[L^1(\lambda)]{T[\La(m,f)]}$. 
Fix $b\in L^\infty(\lambda)$. By Fremlin extension of
Lusin's Theorem \cite[Theorem 2.b]{fremlin_81} (see
also \cite[7.1.13]{bogachev_II}) there exists 
a sequence $\seqn b$ in $\Cts S$ that converges to $b$ in 
$L^1(\lambda)$. However, by property \iref{onto C}, for each 
$n\in\N$ there exists $h_n\in\La_b(m,f)$ such that $T(h_n)=b_n$.

\iref{sw}.
Let $n\in\N$, $\alpha\in\Cts[b]{\R^n}$ 
and $H_1,\ldots,H_n\in\La(m,f)_+$. Consider the
following a.s. equalities
\begin{align*}
\alpha\big(T(H_1),\ldots,T(H_n)\big)
	&=
\lim_k\alpha\big(U(H_1\wedge k),\ldots,U(H_n\wedge k)\big)
	\\&=
\lim_kU\big(\alpha(H_1\wedge k,\ldots,H_n\wedge k)\big)
	\\&=
U\big(\alpha(H_1,\ldots,H_n)\big)
	\\&=
T\big(\alpha(H_1,\ldots,H_n)\big).
\end{align*}
The first one follows from continuity of $\alpha$; the 
second is certainly true for fixed $k$ when $\alpha$ is 
a polynomial and follows from uniform approximation 
by polynomials on the cube $[0,k]^n$; the third one is 
a consequence of the fact that 
$\alpha(H_1\wedge k,\ldots,H_n\wedge k)\in\La_b(m,f)$
and converges to $\alpha(H_1,\ldots,H_n)$ in $L^1(m)$;
the last one is a consequence of $\alpha$ being bounded.

\iref{range}.
Fix $A\in\A(m)$ and let 
\begin{equation}
\Neg_A
	=
\Big\{h\in\La(m,f):
1\ge h\ge 0,\ h(f) \text{ is $m$-null on }A\Big\}.
\end{equation}
It is easy to see that $h\in\Neg_A$ is equivalent to
\begin{align}
\label{null}
0
	=
\int_{B}h(f)dm
	=
\int_{\chi(B)}T(h)d\lambda,
\qquad
B\in\A,\ B\subset A.
\end{align}
By countable additivity the equality 
$0
	=
\int_{\chi(B)}T(h)d\lambda$
extends from $\chi[\A]$ to $\sigma(\clop(S))$ measurable
subsets of $\chi(A)$, and given that, by \iref{meas},  
$T(h)$ is $\sigma(\clop(S))$ measurable, we conclude 
that $T(h)=0$ a.s. on $\chi(A)$. On the other hand this 
last property and \eqref{null} imply $h\in\Neg_A$. 

The set  
$\chi(A)\cap\bigcup_{h\in\Neg_A}\{T(h)>0\}$
is open. Letting $a$ run over the directed family of all 
finite subsets of $\Neg_A$, it follows from $\tau$ 
additivity of $\lambda$ that,
\begin{equation*}
\lambda\Big(\chi(A)\cap\bigcup_{h\in\Neg_A}\{T(h)>0\}\Big)
	=
\lim_a\lambda\Big(\chi(A)\cap\bigcup_{h\in a}\{T(h)>0\}\Big)
	=
0.
\end{equation*}
\end{proof}

By Lemma \ref{lemma T}.\iref{range} the distribution of $h(f)$ under 
$m$ and of $T(h)$ under $\lambda$ coincide save, possibly, 
on a countable set of points. In fact $h(f)\le t$ is equivalent 
to $(h(f)-t)^+=0$ so that, if $h\in\La(m,f)$ and if 
$\{h(f)\le t\}\in\A(m)$, then
\begin{equation}
\label{dist=}
m(h(f)\le t)
	=
\lambda(T(h)\le t).
\end{equation}
Moreover, $\{h(f)\le t\}\in\A(m)$ for all but 
countably many values of $t$.

Given its repeated use, we shall refer to the triple 
$(\chi,\lambda,T)$ satisfying the conditions of Lemma 
\ref{lemma T} as an {\it abstract standard representation} 
of $(m,f)$. It is implicit in the term {\it standard} that 
$\lambda$ is a regular probability defined on the Borel 
subsets of a compact, Hausdorff space. 

\begin{definition}
\label{def abstract}
The triple $(\chi,\lambda,\tilde f)$ is a standard
representation of $(m,f)$ relatively to 
$\Ha\subset\Fun X\cap\La(m,f)$ if
$(\chi,\lambda,T)$ is an abstract standard representation 
of $(m,f)$ and $T(h)=h(\tilde f)$ for all $h\in\Ha$.
\end{definition}

When $f$ is weakly and norm integrable the associated
standard operator $T$ has additional properties, in particular
its restriction $T_0$ of $T$ to $X^*$. In fact,
$T_0\in\Fun{X^*,L^1(\lambda)}$ is weakly compact and
so is then its adjoint, 
$T_0^*\in\Fun{L^\infty(\lambda),X^{**}}$. The restriction 
of $T_0^*$
to the indicators of sets in $\Bor(S)$ gives rise to the
following object.

\begin{definition}
Let $(\chi,\lambda,T)$ be the abstract standard 
representation of $(m,f)$. The standard vector 
measure induced by $(m,f)$ is the unique $X^{**}$ 
valued, countably additive set function of bounded 
variation $\widetilde F$ defined on $\Bor(S)$ 
and satisfying%
\footnote{
The composition $\widetilde F\circo\chi$ is sometimes
referred to as Dunford definite integral of $f$ at $A$.
}
\begin{equation}
\label{tilde F}
\widetilde F\big(\chi(A)\big)x^*
	=
\int_Ax^*fdm,
\qquad
A\in\A.
\end{equation}
\end{definition}

A number of useful properties may be proved easily%
\footnote{
A partial analogue of the equivalence of Pettis integrability
with \iref{w*-to-w} was proved by Edgar 
\cite[Proposition 4.1]{edgar_79}.
}.

\begin{lemma}
\label{lemma stone}
Let $f$ be weakly and norm integrable with abstract 
standard representation $(\chi,\lambda,T)$. Let $T_0$
be the restriction of $T$ to $X^*$ and $\widetilde F$ 
the associated standard vector measure. Then $f$ is
Pettis integrable if and only if either one of the
following equivalent conditions holds:
\begin{enumerate}[(i).]
\item\label{T0 range}
the range of $T_0^*$ belongs to the natural embedding 
of $X$ into $X^{**}$;
\item\label{F range}
$\widetilde F$ takes its values in the natural embedding 
of $X$ into $X^{**}$;
\item\label{w*-to-w}
$T_0$ is weak$^*$-to-weak continuous;
\item\label{b w*-to-w}
$T_0$ is bounded weak$^*$-to-weak continuous%
\footnote{
See \cite[V.5.3]{bible} for a definition of the bounded 
weak$^*$ topology of $X^*$.
}.
\end{enumerate}
\end{lemma}

\begin{proof}
Assume that $f$ is Pettis integrable. The set function 
$F=\widetilde F\circo\chi$ on $\A$ takes its values in
$X$. By linearity $F$ may be extended to an $X$ valued
linear operator defined on the class of all $\A$ simple
functions (endowed with the supremum norm) and such 
extension has its norm dominated by $\int\norm fdm$. 
This permits a further extension, $F_1$, to the closure 
$\B(\A)$ of $\A$ simple functions in the topology of 
uniform convergence. Fix $b\in L^\infty(\lambda)$. By 
Lemma \ref{lemma T}.\iref{dense} there exists a sequence 
$\seqn h$ in $\La_b(m,f)$ such that $T(h_n)$ converges 
to $b$ in $L^1(\lambda)$ and, with no loss of generality, 
such that $T(h_n)$ is bounded by $\norm[L^\infty(\lambda)]{b}$.
But then
\begin{align*}
\bnorm{F_1(h_n(f))-F_1(h_{n+k}(f))}
	&\le
\sup_{x^*\in\mathbb B^*}
\int\abs{x^*f}\babs{h_n(f)-h_{n+k}(f)}dm
	\\&\le
\int T(\norm \cdot)\ \babs{T(h_n-h_{n+k})}d\lambda.
\end{align*}
The sequence $F_1(h_n(f))$ is thus Cauchy in $X$ and 
\begin{align*}
x^*\lim_nF_1(h_n(f))
=
\lim_nx^*F_1(h_n(f))
=
\int bT(x^*)d\lambda
=
T_0^*(b)x^*.
\end{align*}
Thus \iref{T0 range} holds. The implication 
\imply{T0 range}{F range}
is obvious while its converse follows from continuity of $T^*$
and density of simple $\Bor(S)$ measurable functions into
$L^\infty(\lambda)$. Assume \iref{T0 range}. Then we get
\begin{equation*}
x^*x(b)
=
\int bT(x^*)d\lambda,
\qquad
b\in L^\infty(\lambda),\ x^*\in X^*
\end{equation*}
where $x(b)$ is the element of $X$ corresponding to 
$T_0^*(b)$ through the natural embedding of $X$ 
into $X^{**}$. From this it is immediate that $T_0$ 
is weak$^*$-to-weak continuous and that $f$ is Pettis 
integrable. The equivalence of \iref{w*-to-w} with 
\iref{b w*-to-w} follows from the definition of weak 
topology and from \cite[V.5.6]{bible}. 
\end{proof}

Because of the equivalence with \iref{b w*-to-w}
and of \cite[V.5.1]{bible}, if $X$ is separable the property of 
weak$^*$-to-weak continuity of $T_0$ may be proved 
solely in terms of bounded sequences. In other words
separable Banach spaces satisfy the condition of Mazur 
\cite[p. 563]{edgar_79}. In this special case, and with the
additional assumption that $m$ is countably additive, the 
conclusion that all weakly and norm integrable (or norm
bounded) functions are Pettis integrable follows easily 
from the Lebesgue dominated convergence. Another 
obvious case in which Pettis integrability is guaranteed 
is when $X$ is reflexive. This is again immediate from 
\iref{T0 range} and does not require countable additivity 
of $m$.

\section{The Radon-Nikod\'ym property.}
\label{sec RNP}
In this section we consider the implications of the 
Radon-Nykod\'ym property ({\it RNP}) on our results. 

\begin{definition}
\label{def RNP}
$X$ possesses ({\it RNP}) whenever for any countably 
additive, probability space $(S,\Sigma,\nu)$ and any 
measure of finite variation $G\in\Fun{\Sigma,X}$ the 
relation $G\ll\mu$ implies the existence of $g\in\Fun{S,X}$ 
Pettis integrable such that $G(E)=\int_Egd\nu$ for 
each $E\in\Sigma$. 
\end{definition}

\begin{theorem}
\label{th pettis}
Let $X$ possess (RNP) and let $f$ be norm and weakly 
integrable. Then, $f$ is Pettis integrable if and only if 
$(m,f)$ admits a standard representation 
$(\chi,\lambda,\tilde f)$ relative to 
$\Cts[u]{(X,weak)}\cap \La(m,f)$ with 
$\tilde f\in L^1_X(\lambda)$.
Moreover, for each $A\in\A$ there exists a closed subset 
$\tilde A\subset\chi(A)$ such that 
$\lambda\big(\chi(A)\setminus \tilde A\big)
	=
0$ 
and
$\tilde f[\tilde A]\subset\cl[weak]{f[A]}$.
\end{theorem}

\begin{proof}
$f$ is Pettis integrable if and only if the standard vector 
measure $\widetilde F$ induced by $(m,f)$ takes value in 
$X$. Given that $\widetilde F\ll\lambda$ and that $X$ has 
({\it RNP}), this is in turn equivalent to the existence
of $\tilde f\in L^1_X(\lambda)$ such that
\begin{equation}
\label{T=f}
x^*\widetilde F(E)
	=
\int_Ex^*\tilde fd\lambda,
\qquad
E\in\Bor(S).
\end{equation}
By \eqref{T commute}, the a.s. equality $T(x^*)=x^*\tilde f$ 
that follows from \eqref{T=f} extends from $X^*$ to 
$\{\alpha(x^*_1,\ldots,x^*_n):
n\in\N,\ 
\alpha\in\Cts[b]{\R^n}\}$ 
so that $(\chi,\lambda,\tilde f)$ is a standard representation 
of $(m,f)$ relatively to $\Cts[u]{(X,weak)}\cap\La(m,f)$. The 
last claim follows immediately from Theorem \ref{th cover} 
once noted that the restriction of $f$ to any $A\in\A$ is
a representation of the restriction of $\tilde f$ to $\chi(A)$
relatively to $\Cts[ub]{(X,weak)}$.
\end{proof}

An easy implication of Theorem \ref{th pettis} is the
equality
\begin{equation}
\lambda(\tilde f\in E)
	=
m(f\in E),
\qquad
E\subset X,\ E\text{ weakly closed },\ 
f^{-1}(E)\in\A
\end{equation}
which extends the remark on Lemma \ref{lemma T} from the
distribution of $h(f)$ to that of $f$, over weakly open
subsets of $X$.

A much more interesting representation obtains under a 
convenient assumption on the range of $f$. 

\begin{theorem}
\label{th pettis choice}
Let $X$ possess (RNP). Let $f$ be norm and Pettis 
integrable and have essentially weakly closed range%
\footnote{
By this we mean that for each $n\in\N$ there exists 
$A_n\in\A$ with $m(A_n^c)<2^{-n}$ such that 
$f[A_n]$ is weakly closed.
}. 
Define
$\Ha=\Cts[u]{(X,weak)}\cap\La(m,f)$.
There exists $\mu\in\Prob[ca]{f^{-1}[\Bor(X)]}$ such that
\begin{equation}
h(f)\in L^1(\mu)
\qand
\int h(f)dm
	=
\int h(f)d\mu,
\qquad
h\in\Ha.
\end{equation}
\end{theorem}

\begin{proof}
Let $(\chi,\lambda,\tilde f)$ be the standard representation 
of $(m,f)$ relative to $\Ha$, established in Theorem 
\ref{th pettis}. For each $n\in\N$, let $A_n$ be as in the claim
and let $\tilde A_n$ be related with $A_n$ in the same way 
as is $\tilde A$ with $A$ in Theorem \ref{th pettis}. Then, 
$\tilde f\big[\bigcup_n\tilde A_n\big]
\subset 
f\big[\bigcup_nA_n\big]$.
Fix $\omega_0\in\Omega$ and let $\tau(s)=\omega_0$ 
for each $s\notin\bigcup_n\tilde A_n$ or else let
$\tau(s)\in \{f=\tilde f(s)\}$ when $s\in\bigcup_n\tilde A_n$. 
It is then clear that
\begin{align*}
\lambda^*\big(\tilde f\ne f\circo \tau\big)
	\le
\lambda\Big(\bigcap_n\tilde A_n^c\Big)
	\le
\lambda(\tilde A_n^c)
	=
\lambda(\chi(A_n)^c)
	=
m(A_n^c)
	\le
2^{-n}.
\end{align*}
Thus $f\circo\tau\in L^1_X(\lambda)$ and, necessarily, 
$\tau^{-1}(E)\in\Bor(S)(\lambda)$ for every
$E\in f^{-1}[\Bor(X)]$. It follows that
\begin{equation}
\mu
	=
\lambda\circo\tau^{-1}
	\in
\bProb[ca]{f^{-1}[\Bor(X)]}.
\end{equation}
Notice that $h(f)$ is $(f^{-1}[\Bor(X)],\Bor(\R))$ measurable
when $h\in\Cts[u]{(X,weak)}$. The claim then follows from 
an application of the change of variable formula.
\end{proof}

\section{The Pettis integrability property.}
\label{sec PIP}

We return on this property which will be useful in our 
version of Choquet Theorem.

\begin{definition}
Given a finitely additive probability space $(S,\Sigma,\nu)$, 
a Banach space $X$ has ($\nu$-{\it PIP}) if every weakly 
and norm $\nu$-integrable function is Pettis $\nu$ integrable 
too. A Banach space satisfies ({\it PIP}) if the preceding 
property holds for any finitely additive probability space.
\end{definition}

This property, which we briefly considered after Lemma 
\ref{lemma stone}, has been  discussed extensively in the 
literature, in particular by Edgar \cite{edgar_77} and 
\cite{edgar_79} and by Fremlin and Talagrand 
\cite{fremlin_talagrand} (but see also \cite{stefansson} 
and \cite{musial_2002}), although just for the countably 
additive case in which it is, intuitively, a less restrictive 
property. 

Apart from the obvious cases mentioned after Lemma 
\ref{lemma stone}, Edgar \cite[Proposition 3.1]{edgar_79} 
proved that if $X$ is $\nu$-measure compact (for example 
if the weak topology of $X$ is Lindel\"of) then it satisfies the
($\nu$-{\it PIP}). A result in the negative was obtained by 
Fremlin and Talagrand \cite[Thorem 2B]{fremlin_talagrand} 
who proved that $X=\ell^\infty$ fails to possess the countably
additive ({\it PIP})%
\footnote{
Edgar \cite{edgar_77} studies conditions under which 
weakly measurable functions are weakly equivalent to 
measurable ones.
}.

The following result focuses on separability of $X$ 
and establishes only a partial analogue of the countably 
additive case. Nevertheless it permits to decompose
$T$ into a Pettis representable part and a purely non 
representable part. Of course by Theorem 
\ref{th bochner} separability may be replaced with the 
assumption that $f$ is measurable.

\begin{theorem}
\label{th sep}
Let $X$ be separable and $f$ norm 
integrable with abstract, standard representation 
$(\chi,\lambda,T)$. There exist $\eta\in L^1(\lambda)_+$ 
and $\tilde f\in L^1_X(\lambda)$ such that
\begin{equation}
\label{T decomposition}
T(h)
=
h(\tilde f)+T^\perp(h),
\qquad
h\in \La(m,f)
\end{equation}
where $\tilde f=0$ on $\{\eta>0\}$ and
\begin{equation}
T^\perp(x^*)
\ge
\eta
+
\inf_{\sigma\in\Seq{x^*}}\limsup_n T^\perp(\sigma_n),
\qquad
x^*\in X^*
\end{equation}
with $\Seq{x^*}$ denoting bounded sequences in 
$X^*$ which converge weakly$^*$ to $x^*$.
\end{theorem}

\begin{proof}
Since $X$ is separable, $\mathbb B^*$ is metrizable 
in the weak$^*$ topology and $X^*$ admits a 
countable, rational vector space $X_0^*$ which is 
weakly$^*$ dense in $X^*$. Let 
$\mathbb B_0^*
	=
\mathbb B^*\cap X^*_0$.
Select a subset $S_0\subset S$ of full $\lambda$ 
measure such that the functionals $T_s$ are linear 
on $X_0^*$ and that
$\sup_{x^*\in\mathbb B^*_0}\abs{T_s(x^*)}
	\le 
T_s(\norm\cdot)
	<
+\infty$ 
for each $s\in S_0$.

Define $\eta\in\Fun{S}$ by setting $\eta_s=0$
if $s\notin S_0$ or else
\begin{equation}
\eta_s
	=
\sup_\sigma\inf_nT_s(\sigma_n),
\qquad
s\in S_0
\end{equation}
where the $\sup$ is computed over all sequences
$\sigma=\seqn\sigma$ in $\mathbb B^*_0$ which
converge weakly$^*$ to $0$. Clearly, $\eta\ge0$. 

Fix $a\in\R$. If $\seq nk$ is a sequence in $\N$, define 
\begin{equation}
\label{souslin}
A^a_{n_1,\ldots,n_k}
	=
\bigcap_{i=1}^k\{T(x^*_{n_i})\ge a\}
\qtext{if}
d(x^*_{n_i},0)\le2^{-i},
\quad
i=1,\ldots,k
\end{equation}
or else 
$A^a_{n_1,\ldots,n_k}
	=
\emp$. 
The inequality $T(x^*_{n_i})\ge a$ is unaltered if we 
replace $T(x^*_{n_i})$ with its truncation 
$T(x^*_{n_i})\wedge(a+1)\vee(a-1)$. However, by 
Lemma \ref{lemma T}, such truncation is just the 
image under $T$ of 
$x^*_{n_i}\wedge(a+1)\vee(a-1)
	\in
\La_b(m,f)$ 
and is thus continuous. Thus the sets $A^a_{n_1,\ldots,n_k}$ 
are actually closed and define a Souslin scheme so that
\begin{equation}
A^a
	=
\bigcup_{\seq nk}\bigcap_{i=1}^\infty A^a_{n_1,\ldots,n_i}
\end{equation}
is a Souslin set. Moreover,
\begin{equation}
\{\eta\ge a\}\cap S_0
=
\bigcap_{n}A^{a-2^{-n}}\cap S_0
\end{equation}
so that $\eta$ is universally measurable. 

It is easily seen that for every $s\in S_0\cap\{\eta=0\}$ 
$\lim_nT_s(x^*_n)$ exists for every bounded sequence 
$\seqn{x^*}$ in $X_0^*$ which weakly$^*$ converges 
to some $x^*\in X^*$ and that the limit is independent 
of the intervening sequence. We may thus define 
\begin{equation}
q_s(x^*)
=
\lim_nT_s(x^*_n)
\qquad
s\in S_0\cap\{\eta=0\}
\end{equation}
Then, $q_s$ is a linear on $X^*$.

Fix $t>0$ and let $\seqn{z^*}$ be a weakly$^*$ 
convergent sequence in $t\mathbb B^*\cap q_s^{-1}(0)$ 
with limit $z^*\in t\mathbb B^*$. For each $n\in\N$ we 
can find a sequence $\sseq{y^*_{n,k}}{k}$ in
$t\mathbb B_0^*$ weakly$^*$ convergent to
$z^*_n$ and such that $\abs{q_s(y^*_{n,k})}<2^{-n}$.
Using a diagonal argument we can then extract
a sequence $\seq{y^*}{i}$ again in $t\mathbb B^*_0$
which converges weakly$^*$ to $z^*$ and with 
$\abs{q_s(y^*_i)}<2^{-i}$. This implies that
$t\mathbb B^*\cap q_s^{-1}(0)$ is weakly$^*$
closed for each $t>0$ and, by \cite[V.5.6-7]{bible}, 
that $q_s$ is continuous in the weak$^*$ topology. 
But then there exists $\tilde f_s\in X$ such that 
$q_s(x^*)
=
x^*\tilde f_s$
for each $s\in S_0\cap\{\eta=0\}$
and $x^*\in X^*$ and thus such that
\begin{equation}
x^*\tilde f_s
	=
T_s(x^*),
\qquad
x^*\in X^*_0,\ 
s\in S_0\cap \{\eta=0\}.
\end{equation}
Define $\tilde f_s=0$ when $s\notin S_0\cap\{\eta=0\}$.
Clearly, $\tilde f$ is weakly measurable and with separable 
range and thus measurable. Moreover,
$\norm{\tilde f}
	\le
\sup_{x^*\in\mathbb B^*}\abs{T(x^*)}
	\le
T(\norm\cdot)$
so that $\tilde f$ is norm integrable and therefore
$\tilde f\in L^1_X(\lambda)$.

Define $T^\perp(h)=T(h)-h(\tilde f)$. 
If $\seqn{x^*}$ is a bounded sequence in $X_0^*$
weakly$^*$ converging to $0$ and $x^*\in X^*$ then,
\begin{align*}
T^\perp(x^*)
	&\ge 
\inf_kT^\perp(x^*_k)+T^\perp(x^*-x^*_n)
	\\&=
\sset{\eta>0}\inf_kT(x^*_k)+\limsup_nT^\perp(x^*-x^*_n)
	\\&\ge
\eta+\inf_{\sigma\in\Seq{x^*}}\limsup_nT^\perp(\sigma_n).
\end{align*}
\end{proof}

\section{A non compact and non convex 
Choquet theorem}
\label{sec choquet}

To close, we apply our techniques to Choquet integral 
representation on non compact sets. The main result 
in this direction was proved by Edgar \cite{edgar} 
under the assumption that the Banach space $X$ is 
separable and satisfies the ({\it RNP}) (and, as a 
consequence, the Krein-Milman property ({\it KMP}) 
and, as noted above, the countably additive ({\it PIP})). 
In addition, it is assumed that the intervening set
is closed, bounded and convex. One key point in
Choquet Theorem is proving that the set of extreme
points is measurable in some appropriate sense. In
our approach we shall not assume convexity and at
the same time we may disregard measurability issues
to some extent.

\begin{proposition}
\label{pro choquet}
Let $X$ be a Banach space satisfying  (PIP) and (RNP) 
and $C$ a closed and bounded subset of $X$. Fix 
$\Phi\subset\Cts{(X,weak)}$. Assume that $D\subset C$ 
satisfies either one of the following properties:
\begin{enumerate}[(i).]
\item
$(D,weak)$ is a normal topological space%
\footnote{
This condition is satisfied if $X$ is weakly Lindel\"of. 
More examples are contained in the classical work of 
Corson \cite{corson}.
}
or
\item
each $z\in D$ is a point of continuity of $C$%
\footnote{
For example, if $X$ has Kadec norm and $D$ is the intersection
of $C$ with the unit sphere.
}.
\end{enumerate}
Then the condition
\begin{equation}
\label{D dom}
\abs{\phi(x)}
	\le
\sup_{z\in D}\abs{\phi(z)}
	<
+\infty,
\qquad
x\in C,\ \phi\in\Phi
\end{equation}
is satisfied if and only if for each $x\in C$ there exists a 
unique  $\mu_x\in\Prob[ca]{\Bor(D)}$ such that 
\begin{equation}
\label{edgar}
\phi(x)
	=
\int_D \phi(z)d\mu_x,
\qquad
\phi\in\Phi.
\end{equation}
\end{proposition}

\begin{proof}
Of course \eqref{edgar} implies \eqref{D dom}. Assume 
\eqref{D dom}. Then each $x\in C$ acts as a continuous 
linear functional on the linear subspace of $\Cts[b]D$ 
spanned by the elements of $\Phi$. The norm preserving 
extension of each such functional to the whole of 
$\Cts[b]D$ may be represented via a unique regular 
probability $m_x\in\Prob{\Bor(D)}$, \cite[IV.6.2]{bible}. 
In other words we obtain
\begin{equation}
\label{aff}
\phi(f)\in L^1(m_x)
\qand
\phi(x)
	=
\int_D \phi(f)dm_x,
\qquad
x\in C,\ 
\phi\in\Phi
\end{equation}
where $f\in\Fun{D,X}$, the identity map, is weakly 
integrable with respect to $m_x$ (because 
$\mathbb B^*\subset\Cts[b]D$) and bounded in norm. 
Therefore, by ({\it PIP}), $f$ is Pettis integrable as well. 
Under \tiref{i} the restriction of $m_x$ to the Borel 
$\sigma$ algebra $\Bor((D,weak))$ generated by the weakly 
open sets is itself regular so that $f$ has approximately
weakly closed range. On the other hand, if $D$ consists
of the points of continuity of $C$,  then the weak and the
norm topology coincide in restriction to $D$ and again
$f$ has approximately weakly closed range. Moreover, 
$f^{-1}[\Bor(X)]=\Bor(D)$. By Theorem \ref{th pettis choice} 
there exists $\mu_x\in\Prob[ca]{\Bor(D)}$ such that for
every $h\in\Cts[u]{(X,weak)}$
\begin{equation}
h(f)\in L^1(\mu_x)
\qand
\int h(f)dm_x
	=
\int h(f)d\mu_x
\end{equation}
which implies \eqref{edgar}. Given that 
$\{h(f):h\in\Cts[u]{(X,weak)}\}$
is a lattice containing the constant and that the $\sigma$ 
algebra induced by it coincides with $f^{-1}[\Bor(X)]$, 
uniqueness of $\mu_x$ follows from Daniell Theorem.
\end{proof}

By comparison with the original result proved by Edgar, 
Proposition \ref{pro choquet} may appear significantly 
more general but its generality is hindered by the need 
to explicitly assume  ({\it PIP}). Notice that if $C$ is convex 
then $D$ may be the set of its extremal points or of its 
denting points. The existence of such sets is guaranteed 
since, as is well known, the ({\it RNP}) implies the 
({\it KMP}) as well as dentability of bounded sets, 
\cite[Corollary]{davis_phelps}.  In particular if $C$ is 
bounded, closed and convex, and $D$ is the set of its 
denting points, then property \tiref{ii} is satisfied, see 
\cite[Theorem]{lin_et_al}.

\BIB{abbrv}

\end{document}